\documentclass[reqno, 10pt]{amsart}
\theoremstyle{definition}
\newtheorem{definition}{Definition}[section]

\theoremstyle{plain}
\newtheorem{theorem}[definition]{Theorem}

\newtheorem{proposition}[definition]{Proposition}

\theoremstyle{remark}

\newtheorem{remark}[definition]{Remark}

\def\to{\longrightarrow}
\def\comdeg#1{deg(#1)}   
\def\codlev#1{lev(#1)}   

\begin{document}


\title[Combinatorial Aspects of Code Loops]
{Combinatorial Aspects of Code Loops}
\author{Petr Vojt\v echovsk\'y}
\address{Department of Mathematics, Iowa State University, Ames, IA 50011, U.S.A.}
\email{petr@iastate.edu}

\begin{abstract}
The existence and uniqueness (up to equivalence defined below) of code loops
was first established by R.~Griess in \cite{Griess}. Nevertheless, the explicit
construction of code loops remained open until T.Hsu introduced the notion of
symplectic cubic spaces and their Frattini extensions, and pointed out how the
construction of code loops followed from the (purely combinatorial) result of
O.~Chein and E.~Goodaire contained in \cite{CG}. Within this paper, we focus on
their combinatorial construction and prove a more general result
\ref{MainTheorem} using the language of derived forms.
\end{abstract}

\maketitle


\section{Symplectic Cubic Spaces and Code Loops}

\noindent Throughout this paper, let $F=\{0$, $1\}$ be the two-element field,
and let $V$ be a finite-dimensional vector space over $F$. For $v\in V$, let
$|v|$ denote the number of non-zero coordinates of $v$---the {\it weight} of
$v$. When $w$ is another vector in $V$, let $v*w$ denote the vector whose $i$th
coordinate is non-zero if and only if the $i$th coordinate of both $v$ and $w$
is non-zero. A binary linear code $C\le V$ is said to be of {\it level} $r$ if
$r$ is the biggest integer such that $2^r$ divides the weight of every codeword
of $C$. We write $\codlev{C}=r$. A code $C$ is {\it doubly even} if
$\codlev{C}\ge 2$.

For the rest of this section, let $C$ be a doubly even code. Following Griess,
a mapping $\varphi:C\times C\to F$ is called a {\it factor set} if $\varphi(c$,
$c)=|c|/4$, $\varphi(c$, $d)+\varphi(d$, $c)=|c*d|/2$, and $\varphi(c$,
$d)+\varphi(c$, $d+e)+\varphi(d$, $e)+\varphi(c+d$, $e)=|c*d*e|$ is satisfied
for all $c$, $d$, $e\in C$. When $\varphi$ is a factor set, then $(F\times C$,
$\circ)$ with multiplication
\begin{displaymath}
    (\alpha, c)\circ(\beta, d)=(\alpha+\beta+\varphi(c, d), c+d)
\end{displaymath}
becomes a Moufang loop, a {\it code loop} of $C$. R.~Griess shows in
\cite{Griess} that every $C$ admits a factor set $\varphi$, and thus that there
is a code loop for every doubly even code $C$. Moreover, when $\varphi$, $\psi$
are two factor sets for $C$, then they are {\it equivalent} in the sense that
the second derived form $(\varphi+\psi)_2$ is the zero mapping. (See section
$2$ for the definition of derived forms.)

Note that a loop $L$ is a code loop of $C$ if there is a two-element central
subgroup $Z\le Z(L)$ such that $L/Z$ is isomorphic to $C$ as an elementary
abelian $2$-group.

The following ideas are due to T.~Hsu \cite{Hsu}. Let $L$ be a code loop of
$C$. Let $[\gamma$, $\delta]$ denote the commutator of $\gamma$, $\delta$, and
$[\gamma$, $\delta$, $\epsilon]$ the associator of $\gamma$, $\delta$,
$\epsilon\in L$. Define functions $\sigma:C\to Z$, $\chi:C\times C\to Z$, and
$\alpha:C\times C\times C\to Z$ by
\begin{eqnarray}
    \sigma(c)&=&\gamma^2,\nonumber\\
    \chi(c, d)&=&[\gamma, \delta],\label{doublestar}\\
    \alpha(c, d, e)&=&[\gamma, \delta, \epsilon],\nonumber
\end{eqnarray}
where $\gamma$, $\delta$, $\epsilon\in L$ are any preimages of $c$, $d$, $e$
with respect to $L\to L/Z=C$, respectively. One can check that these functions
are well defined, and that the following equalities are satisfied for any $c$,
$d$, $e$, $f\in C$, $n\in \mathbb{N}$ (c.f. Theorems $3.3$ and $4.6$ of
\cite{Hsu}):
\begin{eqnarray}
    \sigma(nc)&=&n\sigma(c), \nonumber\\
    \sigma(c+d)&=&\sigma(c)+\sigma(d)+\chi(c, d),\nonumber\\
    \chi(c, c)&=&0, \nonumber\\
    \chi(c, d)&=&-\chi(d,c),\nonumber\\
    \chi(nc,d)&=&n\chi(c,d), \label{singlestar}\\
    \chi(c+d, e)&=&\chi(c,e)+\chi(d,e)+\alpha(c, d, e),\nonumber\\
    \alpha(c, d, d)&=&\alpha(d, c, c)=\alpha(d, d, c)=0, \nonumber\\
    \alpha(c, d, e)&=&-\alpha(d, c, e)=\alpha(d, e, c), \nonumber\\
    \alpha(nc, d, e)&=&n\alpha(c, d, e), \nonumber\\
    \alpha(c+d, e, f)&=&\alpha(c, e, f)+\alpha(d, e, f), \nonumber
\end{eqnarray}
where the operation in $Z$ is written additively.

The above situation is a special instance of a so-called {\it symplectic cubic
space} $(V$, $\sigma$, $\chi$, $\alpha)$, where $V$ is a vector space over $F$,
and $\sigma: V\to \mathbb Z_2$, $\chi:V\times V\to \mathbb Z_2$,
$\alpha:V\times V\times V\to \mathbb Z_2$ are mappings satisfying
(\ref{singlestar}).

For any symplectic cubic space $(V$, $\sigma$, $\chi$, $\alpha)$ it is
reasonable to define a {\it Frattini extension} $L$, which is a loop with
two-element central subgroup $Z$ such that $L/Z$ is isomorphic to $V$, and such
that $\gamma^2=\sigma(c)$, $[\gamma$, $\delta]=\chi(c$, $d)$, and $[\gamma$,
$\delta$, $\epsilon]=\alpha(c$, $d$, $e)$ is satisfied for all $\gamma$,
$\delta$, $\epsilon$ in $L$. The existence and uniqueness of Frattini
extensions is discussed in detail in \cite{Hsu}. For our purposes it is
sufficient to show that the code loops are precisely the Frattini extensions of
symplectic cubic spaces.

To see this, let $L$ be a code loop of $C$. Then $L$ is a Frattini extension of
the symplectic cubic space $(C$, $\sigma$, $\chi$, $\alpha)$, where $\sigma$,
$\chi$, and $\alpha$ are defined as in (\ref{doublestar}).

Conversely, let $L$ be a Frattini extension of $(V$, $\sigma$, $\chi$,
$\alpha)$. As remarked by T.~Hsu; O.~Chein and E.~Goodaire proved in \cite{CG}
that for any symplectic cubic space $(V$, $\sigma$, $\chi$, $\alpha)$ there is
a doubly even code $C$ isomorphic to $V$, and a basis $\{e_1$, $\dots$, $e_n\}$
of $C$ such that $|e_i|/4=\sigma(e_i)$, $|e_i*e_j|/2=\chi(e_i$, $e_j)$, and
$|e_i*e_j*e_k|=\alpha(e_i$, $e_j$, $e_k)$ for all basis elements $e_i$, $e_j$,
$e_k$. All we have to check then is that $\sigma':c\mapsto |c|/4$, $\chi':(c$,
$d)\mapsto |c*d|/2$, and $\alpha': (c$, $d$, $e)\mapsto |c*d*e|$
form---together with $C$--- a symplectic cubic space $(C$, $\sigma'$, $\chi'$,
$\alpha')$, since then $L$ is a Frattini extension of $(C$, $\sigma'$, $\chi'$,
$\alpha')$, too, and whence $L$ is a code loop of $C$. It is straightforward to
show that $\sigma'$, $\chi'$, and $\alpha'$ satisfy (\ref{singlestar}).


\section{Generalization}

\noindent We have seen in the previous section that code loops can be
characterized as Frattini extensions of symplectic cubic spaces. The crucial
step in the proof was to show that any symplectic cubic space can be identified
with $(C$, $\sigma'$, $\chi'$, $\alpha')$, where $C$ is a certain doubly even
code, and $\sigma'$, $\chi'$, and $\alpha'$ are defined as above. We need to
introduce more notation in order to generalize this result.

Let $I=\{v_1$, $\dots$, $v_s\}$ be a subset of $V$ with possible repetitions.
Then $\sum I$ is defined to be the vector $v_1+\dots +v_s$, and $\prod I$
stands for $v_1*\dots*v_s$. To avoid inconvenience, let $\sum
\emptyset=\prod\emptyset=0$, where $\emptyset$ denotes the empty set.

When $P:V\to F$ is a mapping with $P(0)=0$, M.~Aschbacher defined in
\cite{Asch} the {\it $s$th derived form} $P_s:V^s\to F$ of $P$ by
\begin{equation*}
P_s(v_1, \dots, v_s)=\sum_{J\subseteq I}P\Bigl(\sum J\Bigr).
\end{equation*}
See \cite{Asch}, Section $11$ for the basic properties of derived forms. At
this point, let us at least recall that the derived forms of $P$ can be defined
inductively by
\begin{equation}\label{IndEq}
    P_{s+1}(u, v, v_2, \dots, v_s)=P_s(u, v_2, \dots, v_s)
        +P_s(v, v_2, \dots,v_s)+P_s(u+v, v_2, \dots, v_s).
\end{equation}
The smallest integer $r$ such that $P_s$ is the zero mapping for all $s>r$ is
called the {\it combinatorial degree} of $P$, $\comdeg{P}$. Such an integer is
guaranteed to exist and cannot exceed the dimension of $V$.

Since $\sigma'$, $\chi'$, and $\alpha'$ are related by
polarization---$\sigma'(c+d)=\sigma'(c)+\sigma'(d)+\chi'(c$, $d)$, $\chi'(c+d$,
$e)=\chi'(c$, $e)+\chi'(d$, $e) +\alpha'(c$, $d$, $e)$---we see that
$\chi'=\sigma'_2$, and $\alpha'=\sigma'_3$. Therefore the Chein's and
~Goodaire's result can be restated as follows:

{\it Given $P:V\to F$ with $P(0)=0$, $\comdeg{P}=3$, there is a doubly even
code $C$ isomorphic to $V$ such that $P(c)=|c|/4$ for all $c\in C$.}

In the rest of the paper, we prove the main result:

\begin{theorem}\label{MainTheorem}
Let $P:V\to F$ be a mapping of combinatorial degree $r+1$. Then there is
a binary linear code $C$ of level $r$ isomorphic to $V$ such that $P(c)=|c|/2^r$
is satisfied for each codeword $c$ in $C$.
\end{theorem}


\section{Constructing Binary Linear Codes From derived Forms}

\noindent For the sake of brevity let us write $P(I)$ instead of $P_s(v_1$,
$\dots$, $v_s)$, where still $I=\{v_1$, $\dots$, $v_s\}$. Let $P(\emptyset)=0$.
Using this notation, the reverse formula for derived forms can be elegantly
written as
\begin{equation}\label{ReverseFormula}
    P\Bigl(\sum I\Bigr)=\sum_{J\subseteq I} P(J).
\end{equation}
This follows from (\ref{IndEq}) by induction on $|I|$.

Also recall the explicit formulae for the weights of sums and products of
vectors in V:
\begin{eqnarray}
    \Bigl|\sum I\Bigr|=\sum_{J\subseteq I}(-2)^{|J|-1}\Bigl|\prod J\Bigr|, \label{SumToProduct}\\
    2^{s-1}\Bigl|\prod I\Bigr|=\sum_{J\subseteq I}(-1)^{|J|-1}\Bigl|\sum J\Bigr|. \label{ProductToSum}
\end{eqnarray}

\begin{proposition}\label{Equivalence}
Let $P:V\to F$ be a mapping with $P(0)=0$. The following conditions are equivalent:
\begin{enumerate}
\item[(i)] $2^rP(\sum I)\equiv |\sum I|\pmod{2^{r+1}}$ for any subset $I\subseteq
V$,
\item[(ii)] $2^{r-|I|+1}P(I)\equiv |\prod I|\pmod{2^{r-|I|+2}}$
for any subset $I\subseteq V$.
\end{enumerate}
\begin{proof}
Suppose $(i)$ is satisfied. Let $I$ be a subset of $V$.
We have
\begin{equation*}
    P(I)\equiv\sum_{J\subseteq I}P\Big(\sum J\Big)\equiv \sum_{J\subseteq I} (-1)^{|J|-1} P\Big(\sum J\Big)\pmod2.
\end{equation*}
Multiplying this congruence by $2^r$, and using $(i)$, we immediately obtain
\begin{equation*}
    2^rP(I)\equiv \sum_{J\subseteq I} (-1)^{|J|-1} \Big|\sum J\Big|\pmod{2^{r+1}}.
\end{equation*}
Using (\ref{ProductToSum}), we finally get
\begin{equation*}
    2^{r-|I|+1}P(I)\equiv 2^{1-|I|}\cdot\sum_{J\subseteq I}
    (-1)^{|J|-1} \Big|\sum J\Big| \equiv \Big|\prod I\Big|\pmod{2^{r-|I|+2}}.
\end{equation*}

Now assume that $(ii)$ is satisfied. By the reverse formula
(\ref{ReverseFormula}), and after some convenient rearrangements, we see that
\begin{equation*}
    2^rP\Big(\sum I\Big)\equiv \sum_{J\subseteq I} (-1)^{|J|-1}2^r P(J)\pmod{2^{r+1}}.
\end{equation*}
Condition $(ii)$ says that $2^rP(J)\equiv 2^{|J|-1}|\prod J|\pmod{2^{r+1}}$.
Thanks to (\ref{SumToProduct}), we get
\begin{equation*}
    2^rP\Big(\sum I\Big)\equiv \sum_{J\subseteq I} (-2)^{|J|-1}\Big|\prod J\Big| \equiv \Big|\sum I\Big|\pmod{2^{r+1}},
\end{equation*}
as desired.
\end{proof}
\end{proposition}

Let us first outline the construction of $C$ in words.

Let $\{v_1$, \dots, $v_m\}$ be a basis for $V$. Suppose that we have found
linearly independent vectors $c_1$, \dots, $c_m$, which generate a linear code
$C$ of level $r$. Let us identify $v_i$ with $c_i$, for $1\le i\le m$. Every
codeword $c\in C$ can be expressed as $\sum I$ for some $I\subseteq\{c_1$,
\dots, $c_m\}$. We would like to have $P(\sum I)\equiv|\sum I|/ 2^r\pmod2$ for
every $I$. According to Proposition \ref{Equivalence}, we only need to
guarantee condition
\begin{equation}
    2^{r-|I|+1}P(I)\equiv \Big|\prod I\Big|\pmod{2^{r-|I|+2}} \tag{$I$}
\end{equation}
for every $I\subseteq\{c_1$, \dots, $c_m\}$.

We construct the vectors $c_1$, \dots, $c_m$ in $2^m-1$ steps. Let us label
these steps by non-empty subsets of $\{1$, \dots, $m\}$, and order them as
follows: if $I=\{i_1$, \dots, $i_k\}$, $J=\{j_1$, \dots, $j_l\}$, where $i_1 >
\cdots > i_k$, $j_1 > \cdots > j_l$, then $I\le J$ if and only if $(i_1$,
\dots, $i_k)\le(j_1$, \dots, $j_l)$ lexicographically.

Vector $c_i$ is introduced in step $\{i\}$. In each step, a certain number of
coordinates is adjoint to each of the previously introduced vectors. If
$i\not\in I$, and if $c_i$ has already been mentioned, we extend $c_i$ by zeros
in step $I$. Let us identify subsets of $\{1$, \dots, $m\}$ with subsets of
$\{c_1$, \dots, $c_m\}$ in the natural way. After each step $I$, we check that
all conditions $(J)$, $J\le I$, are satisfied, and that the previously
introduced vectors generate a linear code of level at least $r$. In fact, after
the construction is finished, we necessarily get $\codlev{C}=r$, otherwise $P$
is the zero mapping.

Moreover, note that when $|I|>1$, then all vectors $c_i\in I$ have already been
introduced. In order to make the construction more transparent, we will
construct the vectors in such a way that $\prod I=0$ is satisfied before step
$I$, for $|I|>1$.

Now, we are ready to begin with the construction.

\vskip 3mm

{\bf Steps $\{i\}$:}

Add $2^{r+1}$ coordinates to all previously introduced vectors. Define a new
vector $c_i$ whose only non-zero coordinates are among the last $2^{r+1}$
coordinates, which consist of $2^r$ ones and $2^r$ zeros if $P(v_i)=1$, and of
$2^{r+1}$ ones if $P(v_i)=0$.

Then $2^rP(v_i)\equiv |c_i|\pmod2$, and condition $(J)$ remains valid for every
$J\le\{i\}$ because $i\not\in J$. All introduced vectors generate a linear code
of level at least $r$.

\vskip 3mm

{\bf Steps $I$, for $|I|>1$:}

We need the following rather general proposition. It is the key to the whole
construction.

\begin{proposition}\label{CKey}
Let $W = F^{2^k}$ be a vector space over $F$. Let $0<l+1<k$. There are linearly
independent vectors $w_0$, \dots, $w_l\in W$ such that for every proper subset
$A$ of $\{w_0$, \dots, $w_l\}$ we have $|\prod A|=2^{k-|A|}$, and
$|w_0*\cdots*w_l|=2^{k-l-2}$.
\begin{proof}
First, we define real vectors $u_0$, \dots, $u_l\in F^{2^l}$, where
$u_i=(u_{i,j})_{j=0}^{2^l-1}$, $0\le i\le l$. Let us identify the number
$j=\sum_{i=0}^{l-1}j_i 2^i$ with the vector $(j_0$, \dots, $j_{l-1}) \in F^l$.
Let $j^\perp$ denote the complement of $j$ in $F^l$. Let $\varphi : F^l \to F$
be a mapping defined by $\varphi(j)\equiv |j^\perp|\pmod2$. For $0\le i<l$,
$0\le j<2^l$, put $u_{i,j}=j_i$. For $0\le j<2^l$, define
$u_{l,j}=1/4+1/2\cdot\varphi(j)$, i.e. $u_{l,j}\in \{1/4$, $3/4\}$.

To construct vectors $w_i$ from $u_i$, $0\le i\le l$, replace each $u_{i,j}$
with a block of $2^{k-l}\cdot u_{i,j}$ ones and $2^{k-l}\cdot(1-u_{i,j})$
zeros.

We need to check that vectors $w_0$, \dots, $w_l$ have the desired properties.
Let us get started with $|w_0*\cdots*w_l|$. There is only one coordinate $j$,
namely $2^l-1$, for which $1=j_i=u_{i,j}$, $0\le i< l$. Since $\varphi(j)=0$,
we have $u_{l,j}=1/4$. Therefore $|w_0*\cdots*w_l|=2^{k-l-2}$.

Let $A$ be a proper subset of $\{w_0$, \dots, $w_l\}$. Suppose, for a while,
that $w_l\not\in A$. Define $M=\{0\le j<2^l|$ $u_{i,j}=1$ for all $w_i\in A\}$.
Clearly, $|\prod A|=2^{k-l}|M|$. Because $u_{i,j}$ is arbitrary for $w_i\not\in
A$, we have $|M|=2^{l-|A|}$. In other words, $|\prod A|=2^{k-|A|}$.

Suppose that $w_l\in A$. For $t=0$, $1$, put $M_t=\{0\le j< 2^l|$ $u_{i,j}=1$
for $w_i\in A\setminus \{w_l\}$, and $\varphi(j)=t\}$. Then $M_0\cap
M_1=\emptyset$, and $|M_0\cup M_1|=2^{l-|A|}$. Since $|M_0|=|M_1|=2^{l-1-|A|}$,
we get $|\prod A|=1/4\cdot 2^{k-l}\cdot|M_0|+3/4\cdot 2^{k-l}\cdot|M_1|=
2^{k-l}\cdot 2^{l-|A|}=2^{k-|A|}$.
\end{proof}
\end{proposition}

If $P(I)=0$, we do not need to make any changes. Condition $(I)$ is satisfied
because $|\prod I|=0$ has been true before step $I$.

Suppose that $P(I)=1$. Then $r+1=\comdeg{P}\ge |I|$, and we may use Proposition
\ref{CKey} with parameters $l=|I|-1$, $k=r+2$ to obtain vectors $w_0$, \dots,
$w_{|I|-1}$. We extend vectors from $I$ by these vectors $w_i$, one by one (in
any order). By Proposition \ref{CKey} we have $|\prod
I|=2^{r+2-(|I|-1)-2}=2^{r-|I|+1}P(I)$. Let $J<I$. If $J$ is not a proper subset
of $I$, then $|\prod J|$ did not change (vectors not involved in $I$ are
extended by zeros), and that's why condition $(J)$ still holds. If $J$ is a
proper subset of $I$, then $|\prod J|$ increased by $2^{r+2-|J|}$ (according to
Proposition \ref{CKey}), therefore condition $(J)$ holds, too.

All introduced vectors generate a linear code of level at least $r$, and we are
done.

\begin{remark} For the sake of completeness, let us consider the (much easier)
inverse problem of Theorem \ref{MainTheorem}: given a binary linear code $C$ of
level $r$, construct $P:C\to F$ by $P(c)=|c|/2^r$, $c\in C$. Then
$\comdeg{P}\le r+1$. (See \cite{Asch}, Lemma $11.4$, or \cite{PV}.)
\end{remark}

\bibliographystyle{amsalpha}

\end{document}